\input mssymb

\def\qedf#1{}
\def\qed#1{\hfill\rlap{\hskip\rightskip$\smiley_{#1}$}$\phantom{\smiley{#1}}$}
\def\QED#1{\hfill\rlap{\hskip\rightskip$\SMILEY_{#1}$}$\phantom{\smiley{#1}}$}

\newif\ifdvips
\newif\ifdvitwops

\dvipstrue

\hsize=0.9\hsize
\vsize=0.95\vsize
\baselineskip=1.5\baselineskip
\hfuzz 10pt
\parindent=0em

\def\boxit#1{\setbox0=\hbox{~#1\vphantom{ig1}~}%
\setbox1\hbox to \wd0{\hrulefill}%
\wd1=0cm%
\dimen0\ht0\advance\dimen0 2pt
\dimen1\dp0\advance\dimen1 2pt
\leavevmode\hbox{\vrule\raise\dimen0\copy1
			\lower\dimen1\box1
			\box0
		\vrule}}
\font\circle=lcircle10
\setbox0=\hbox{~~~~~}
\setbox1=\hbox to \wd0{\hfill$\scriptstyle\smile$\hfill} 
\setbox2=\hbox to \wd0{\hfill$\cdot\,\cdot$\hfill} 
\setbox3=\hbox to \wd0{\hfill\hskip4.8pt\circle i\hskip-4.8pt\hfill} 
\setbox5=\hbox to  \wd0{\hfill--\hfill}

\wd1=0cm
\wd2=0cm
\wd3=0cm
\wd4=0cm
\wd5=0cm

\newbox\SMILEBOX
\setbox\SMILEBOX \hbox {\lower 0.4ex\copy1
                 \raise 0.3ex\copy2
                 \raise 0.5ex\copy3
                \copy4
                \copy0{}}

\def\SMILEY{\leavevmode\copy\SMILEBOX}

\newbox\smilebox
\setbox\smilebox \hbox {\lower 0.4ex\box5
                 \raise 0.3ex\box2
                 \raise 0.5ex\box3
                \box4
                \box0{}}
\def\smiley{\leavevmode\copy\smilebox}

\def\on{{\restriction}}

\let\frak\bf
\let\fraktur\frak

\def\level{{\rm Level}}
\def\act{{\rm active}}

\def\c{{\frak c}}
\def\r{{\frak r}}

\def\d{{\fraktur d}}

\def\sux{{\rm succ}}
\def\|#1|{\left\Vert #1\right\Vert}

\def\bareta{\bar\eta}


\def\thinks{\models}


\def\pre#1#2{{}^{#1}\!#2}
\def\erp#1^#2{\pre{#2}{#1}}
\def\to{\rightarrow}

\def\cut{\cap}
\def\extend{\mathord{{}^{\frown}}}

\def\fct{\pre{\omega}{\omega}}

\def\finSeq{\erp\omega^{\lless\omega}}

\def\kzerooneseq#1 {\pre{#1}{2}}

\def\<#1>{\hbox{$\langle #1\rangle$}}  

 \def\stronger{\ge}
 \def\weaker{\le}
 
						\def\lessdot{\mathrel{\mathord{<}\!\!%
  \raise 0.8 pt\hbox{$\scriptstyle\circ$}}}  
\let\lesscirc=\lessdot


\def\dom{\mathop{\rm dom}\nolimits}
\def\rng{\mathop{\rm rng}\nolimits}
\def\min{\mathop{\rm min}}

\def\sll{\rm}
\def\stem{{\sll stem}}

\def\split{{\sll split}}


\def\vu{\nu}

\def\[#1]{\lfloor #1 \rfloor}

\def\lless{\mathord{<}}



\newcount\secno
\newcount\theono   

\catcode`@=11
\newwrite\mgfile

\openin\mgfile \jobname.mg
\ifeof\mgfile
\def\slalomdef{0.1} 
\def\coverdef{0.2} 
\def\invDef{0.3} 
\def\ptwise{1.1} 
\def\obvious{1.2} 
\def\obvRemark{1.3} 
\def\uncount{1.4} 
\def\trivZFC{1.5} 
\def\trivZFCrem{1.6} 
\def\moreZFC{1.7} 
\def\techLemma{1.8} 
\def\productCor{1.9} 
\def\const{1.10} 
\def\constants{1.11} 
\def\diverge{1.12} 
\def\power{1.13} 
\def\addition{1.14} 
\def\multiple{1.15} 
\def\nstarDef{2.1} 
\def\normdef{2.2} 
\def\thingydef{2.3} 
\def\fggRem{2.4} 
\def\Notation{2.5} 
\def\singledef{2.6} 
\def\AxiomA{2.7} 
\def\leqsubn{2.8} 
\def\normnotation{2.9} 
\def\cfCCCXXVI{2.10} 
\def\normrem{2.11} 
\def\olddontcover{2.12} 
\def\mainTHEO{3.1} 
\def\onlytwo{3.2} 
\def\Blassanswer{3.3} 
\def\squareex{3.4} 
\def\iterdef{3.5} 
\def\facts{3.6} 
\def\leveldef{3.7} 
\def\levelrem{3.8} 
\def\splitlevel{3.9} 
\def\extendDf{3.10} 
\def\fkDef{3.11} 
\def\densedef{3.12} 
\def\bartrim{3.13} 
\def\decidedef{3.14} 
\def\activefinite{3.15} 
\def\levelsize{3.16} 
\def\kldef{3.17} 
\def\pruneDef{3.18} 
\def\gamedef{3.19} 
\def\strategy{3.20} 
\def\thinning{3.21} 
\def\densefour{3.22} 
\def\propercoro{3.23} 
\def\properremark{3.24} 
\def\smallsupport{3.25} 
\def\supportcor{3.26} 
\def\plusnotation{3.27} 
\def\smalll{3.28} 
\def\densetwo{3.29} 
\def\densesix{3.30} 
\def\geproof{3.31} 
\def\almostall{3.32} 
\def\leproof{3.33} 
\def \Blass{{\fam \slfam \tensl Blass\/}}
\def \vanDouwen{{\fam \slfam \tensl van Douwen\/}}
\def \ComfortNegrepontis{{\fam \slfam \tensl Comfort-Negrepontis}}
\def \Miller{{\fam \slfam \tensl Miller}}
\def \ShelahCCCXXVI{{\fam \slfam \tensl Shelah 326\/}}
\def \ShelahCDXLVIIIa{{\fam \slfam \tensl Shelah 448a}}
\def \Vaughan{{\fam \slfam \tensl Vaughan\/}}
	\else\closein\mgfile\relax\input \jobname.mg\fi
\relax
\openout\mgfile=\jobname.mg

\newif\ifproofmode
\proofmodefalse            

\def\@nofirst#1{}

\def\neusection{\global\advance\secno by 1\relax \global\theono=0\relax}
\def\neuchap{\global\secno=0\relax\global\theono=0\relax}

\neuchap

\def\labelit#1{\global\advance\theono by 1%
             \global\edef#1{%
             \number\secno.\number\theono}%
             \write\mgfile{\@definition{#1}}%
}


\def\prop#1#2:{%
\labelit{#1}%
\smallbreak\noindent%
\@markit{#1}%
{\bf #2:}}

\def\@definition#1{\string\def\string#1{#1}
\expandafter\@nofirst\string\%
(\the\pageno)}

\def\@markit#1{
\ifproofmode\llap{{\tt \expandafter\@nofirst\string#1\ }}\fi%
{\bf #1 \ }
}
 
\def\labelcomment#1{\write\mgfile{\expandafter
		\@nofirst\string\%---#1}} 

\catcode`@=12



%


\def\begindent{\par\smallskip \begingroup \parindent=2cm}
\def\begindent{\par\smallskip\begingroup\advance\parindent by
1cm\advance\rightskip by 1.5 cm plus 1 cm minus 1 cm\relax}
\def\ite#1 {\item{(#1)} }
\def\endent{\smallskip\endgroup}

\headline={\hfill}



%


\def\mybeginsection#1\par{\vskip 0pt plus1cm \penalty -250 \vskip
0pt plus-1cm \bigskip \vskip \parskip \message {#1}%
\leftline {\underbar{\bf #1}}\nobreak \smallskip \noindent}

\footline={\hfill{\tenrm \folio}\hfill{\tenrm
Revised version,         April 1992}}
\newtoks\filename

\def\nextinput #1 #2\par{\neusection
\mybeginsection{#2}\par\filename={#1}
}

\def\dont#1\par{{\let\input\relax #1
\par}}


%


\def\Forces_#1``#2''{\forces_{#1}\hbox{``#2''}}

\def\name#1{\mathchoice%
{\setbox0=\hbox{$\displaystyle #1$}
\setbox1=\vtop{\ialign{##\crcr
$\hfil{\displaystyle #1}\hfil$\crcr\noalign{\kern2pt%
\nointerlineskip}$\hfil\mathord{\displaystyle \sim}%
\hfil$\crcr\noalign{\kern3pt\nointerlineskip}}}%
\setbox2=\hbox{$\displaystyle \sim$}%
\wd1=\wd0\dp1=0cm\ifdim\wd2>\wd1 \wd1=\wd2\else\relax\fi
\ht1=\ht0\relax
\box1}%
{\setbox0=\hbox{$\textstyle #1$}
\setbox1=\vtop{\ialign{##\crcr
$\hfil{\textstyle #1}\hfil$\crcr\noalign{\kern1.2pt%
\nointerlineskip}$\hfil\mathord{\textstyle \sim}%
\hfil$\crcr\noalign{\kern3pt\nointerlineskip}}}%
\setbox2=\hbox{$\textstyle \sim$}%
\wd1=\wd0\dp1=0cm\ifdim\wd2>\wd1 \wd1=\wd2\else\relax\fi
\ht1=\ht0\relax
\box1}{\setbox0=\hbox{$\scriptstyle #1$}
\setbox1=\vtop{\ialign{##\crcr
$\hfil{\scriptstyle #1}\hfil$\crcr\noalign{\kern1pt%
\nointerlineskip}$\hfil\mathord{\scriptstyle \sim}%
\hfil$\crcr\noalign{\kern2.1pt\nointerlineskip}}}%
\setbox2=\hbox{$\scriptstyle \sim$}%
\wd1=\wd0\dp1=0cm\ifdim\wd2>\wd1 \wd1=\wd2\else\relax\fi
\ht1=\ht0\relax
\box1}{\setbox0=\hbox{$\scriptscriptstyle #1$}
\vtop{\ialign{##\crcr
$\hfil{\scriptscriptstyle #1}\hfil$\crcr\noalign{\kern1pt%
\nointerlineskip}$\hfil\mathord{\scriptscriptstyle \sim}%
\hfil$\crcr\noalign{\kern1.5pt\nointerlineskip}}}%
}}

{\catcode`\@=11\gdef\neubox{\alloc@ 4\box \chardef \insc@unt}}
\def\assign#1{\neubox\current
\setbox\current=\hbox{$\name #1$}
\edef\next{\noexpand\edef
\csname n#1\endcsname{{\noexpand\copy\the\current}}}\next}

\def\nofirst#1{}
\def\assignp#1{\neubox\current
\setbox\current=\hbox{$\name #1$}
\edef\next{\noexpand\edef
\csname
n\expandafter\nofirst
\string#1\endcsname{{\noexpand\copy\the\current}}}\next}

\assignp\tau

\assign f
\assign g
\assign h
\assign q
\assign Q
\assign p
\assignp\sigma
\assignp\alpha
\assign r
\assign x
\assign y
\assign k


\def\compile#1#2{\neubox\current
\setbox\current=\hbox{$#2$}
\edef#1{\noexpand\copy\the\current}}



\let\fraktur\frak

\def\c{{\fraktur c}}   

\def\P{{\fraktur P}}






%
%
%
%

\def\forces{\Vdash}

\def\B{{\cal B}}
\def\trim#1^#2{#1^{[#2]}}

\newcount\condno

\def\setupcond#1{\advance\condno by 1 %
\edef\next{\def\csname#1\endcsname{#1}%
\def\csname x#1\endcsname{\the\condno}}%
\next}

  \def\finite	{finite}
  \def\onepair	{onepair}
  \def\kplusone	{kplusone}
  \def\smallevl	{smalllevel}
  \def\tauonk	{tauonk}
  \def\tauofk	{tau(k)}

  \def\ALMOSTALL{almostall}

\def\romanize#1{\edef\next{\romannumeral#1\relax}%
\edef\nnext{\uppercase{\next}}\nnext}

\def\frame#1{\expandafter\boxt\:{#1}}
\def\boxt(#1){\item{\boxit{#1}}}

\def\:#1{(\def\argument{#1}%
 \ifx\finite   \argument \romanize{1}\else
  \ifx\onepair  \argument \romanize{2}\else
   \ifx\kplusone \argument \romanize{3}\else
    \ifx\smallevl \argument \romanize{4}\else
     \ifx\tauonk   \argument \romanize{5}\else
      \ifx\tauofk   \argument \romanize{5}a\else
       \ifx\ALMOSTALL\argument \romanize{6}\else
		\message{argument not recognized:}\show\argument
       \fi
      \fi
     \fi
    \fi
   \fi
  \fi
 \fi
)}

\def\thingy{progressive}

\def\theset#1:#2\\#3\endset{
\setbox0=\hbox{#2}
\big\{#1: \vtop{\parindent=0cm\hsize=\wd0
\rightskip=0cm plus 0.2\wd0 minus 0.2\wd0
#2\break #3 $\big\}$}}

\advance\medskipamount by 0pt plus 2 pt minus 2pt

\def\proof{\smallskip\noindent {\bf Proof: }}

\def\accountant/{{\sl accountant}}
\def\spendthrift/{{\sl spendthrift}}

	\magnification\mag
	\newcount\EINS
	\newcount\ZEHN
	\newcount\HUNDERT
	\newcount\TAUSEND
	\newcount\scratch
	\newskip\totalheight
\ifdvips \input epsf
\else\ifdvitwops  
	\scratch\mag
	\divide\scratch by 10
	\multiply\scratch by 10
	\EINS=\mag
	\advance\EINS by -\scratch
	
	\scratch\mag
	\divide\scratch by 100
	\multiply\scratch by 100
	\ZEHN=\mag
	\advance\ZEHN by -\scratch
	\divide\ZEHN by 10
	
	\scratch\mag
	\divide\scratch by 1000
	\multiply\scratch by 1000
	\HUNDERT=\mag
	\advance\HUNDERT by -\scratch
	\divide\HUNDERT by 100

	\TAUSEND=\mag
	\divide\TAUSEND by 1000
	
	\edef\scalesize{\number\TAUSEND.\number\HUNDERT\number\ZEHN\number\EINS\space}

	\def\epsffile[#1 #2 #3 #4]#5{%
	\totalheight=#4 pt
	\advance\totalheight by -#2 pt
	\vbox to \totalheight{\vskip#4 pt
	\moveleft#1 pt\hbox to #3 pt{%
	\includegraphics{#5}\hskip#3 pt}%
	\vskip-#2 pt}}%
\else 	\def\epsffile[#1 #2 #3 #4]#5{%
	\totalheight=#4 pt
	\advance\totalheight by -#2 pt
	\vbox to \totalheight{\vfill \centerline{Picture missing}\vfill}}
\fi\fi

\proofmodefalse

\hyphenation{spend-thrift}

\font\large=cmr17

^^L

\centerline{\large Many simple cardinal invariants}
\smallskip
\centerline{November 1991}
\medskip

\centerline{\bf Martin Goldstern\footnote{$^1$}{\rm supported by Israeli
Academy of 
Sciences, Basic Research Fund}}
\smallskip
\centerline{Bar Ilan University}
\bigskip
\centerline{\bf Saharon Shelah\footnote{$^2$}{\rm
\baselineskip=0.5\baselineskip
  Publication 448.  Supported partially by Israeli Academy Of
Sciences, Basic Research Fund and by the Edmund Landau Center
for research in Mathematical Analysis (supported by the Minerva
Foundation (Germany))\endgraf}}
\smallskip
\centerline{Hebrew University of Jerusalem}
\bigskip

{\advance\leftskip 0.15\hsize
\advance \rightskip0.15\hsize plus 1cm
Abstract:  For $g<f$ in $\omega^\omega$ we define $\c(f,g)$ be the
least number of uniform trees with $g$-splitting needed to cover a
uniform tree with $f$-splitting.  We show that we can simultaneously
force $\aleph_1$ many different values for different functions
$(f,g)$. In the language of [\Blass]: There may be $\aleph_1$ many
distinct uniform $\bf\Pi^0_1$ characteristics. 

}

\advance\secno by -1

 \nextinput intro {0. Introduction}

[\Blass] defined a classification of certain cardinal invariants of
the continuum, based
on the Borel hierarchy.  For example,  to every $\Pi^0_1$ formula
$\varphi(x,y) = \forall n R(x\on n, y \on n)$ ($R$ recursive)  the cardinal 
$$ \kappa_{\varphi} := \min \{\B \subseteq \fct: \forall x\in \fct
\exists y \in \B : \varphi(x,y) \}$$
is the ``uniform $\Pi^0_1$ characteristic''   associated to~$\varphi$.  

Blass proved  structure theorems on simple cardinal invariants,
e.g., that there is a smallest $\bf\Pi^0_1$ characteristic
(namely, ${\bf Cov}({\cal M})$, the smallest number
of first category sets needed to cover the reals),
and also that the $\Pi^0_2$-characteristics can
behave quite chaotically.  He asked whether the
known uniform $\Pi^0_1$  characteristics ($\c$, 
$\d$, $\r$, ${\bf Cov}({\cal M})$) are the only ones or (since that is very
unlikely) whether there could be  a  reasonable classification of
the uniform $\Pi^0_1$ characteristics --- say, a small  list that
contains all these invariants.  

In this paper we give a strong negative answer to this question:  For
two $\bf \Pi^0_1$ formulas $\varphi_1$, $\varphi_2$ we say that
$\varphi_1$ and $\varphi_2$ define ``potentially nonequal
characteristics'' if $\kappa_{\varphi_1} \not= \kappa_{\varphi_2}$ is
consistent.   We say that $\varphi_1$ and $\varphi_2 $ define ``actually
different characteristics'', if $\kappa_{\varphi_1}\not= \varphi_2$. 

We will find a family of $\bf\Pi^0_1$-formulas indexed by a real
parameter $(f,g)$, and we will show  not only that there is a perfect
set of parameters  which defines 
pairwise potentially nonequal $\bf \Pi^0_1$-characteristics, but we
produce a single universe in which (at least) $\aleph_1$ many
cardinals appear as  $\bf \Pi^0_1$-characteristics.  (In fact it is
also possible to produce a universe where there is a perfect set of
parameters  defining pairwise actually different  $\bf
\Pi^0_1$-characteristics. See [\ShelahCDXLVIIIa]). 

If we want more than countably many cardinals, we  obviously have
to use the boldface 
pointclass.   But the proof also produces many lightface uniform
$\Pi^0_1$ characteristics.

For more information on cardinal invariants, see  [\Blass],
[\vanDouwen], [\Vaughan]. 

\bigskip

\relax From another point of view, this paper is part of the program
of 
finding consistency  techniques for a large continuum, i.e.,
we want $2^{{\aleph_0}} >\aleph_2$ and have many values for cardinal
invariants.    We use a countable support product of 
forcing notions with an  axiom A structure. 

We will use invariants that were implicitly introduced in
[\ShelahCCCXXVI, \S 2], where it was proved that $\c(f,g)$ and
$\c(f',g')$  (see below) may be distinct.

\bigskip\goodbreak

\prop\slalomdef Definition:  If $f\in\fct$, we say that $\bar B =
\<B_k:k\in \omega>$ is an $f$-slalom 
if for all $k$, $|B_k| = f(k)$.   
We  write $h\in \bar B$ for $h\in
\prod_n B_n$, i.e., $\forall n \, h(n)\in B_n$. (See figure 1)  This is a
$\Pi^0_1$-formula in the variables $h$ and~$\bar B$. \qedf\slalomdef

Some authors call the set $\{h: h \in \bar B\} $ a ``belt'', or
``uniform tree''. 

For example,  $\prod\limits_n f(n)$ is an $f$-slalom,
because we identify the number $f(n)$ with the set of  predecessors,
$\{0, \ldots, f(n)-1\}$. 
\midinsert
$$\vcenter{\epsffile[-7 109 314 272]{448-pic3.ps}}$$
\centerline{Figure 1: A slalom}
\endinsert

\prop\coverdef Definition:  Assume $f,g\in \fct$. Assume that $\B$ is
a family of $g$-slaloms, and $\bar A = \< A_k:k\in \omega>$ is an
$f$-slalom.

We say that $\B$ {\bf covers} $\bar A$ iff: 
$$(\star)\qquad\qquad \hbox{for all $s\in \bar A$ there is $\bar  B\in
\B$ such that } \quad 	s\in \bar B \eqno\qedf\coverdef$$

\prop\invDef Definition:  Assume $f, g\in \fct$.  Then we define the
cardinal invariant $\c({f,g})$ to be the 
minimal number of $g$-slaloms needed to cover an $f$-slalom.  \qedf\invDef 

(Clearly this makes sense only if $\forall k \, f(k), g(k) >0$, so we
will assume that from now on.) 

  This is a uniform $\Pi^0_1 $-characteristic.  (Strictly speaking, we
are  not working in $\fct$, but rather in $\pre \omega
{\big([\omega]^{{<}\omega}\big)}  $, but a trivial coding translates
$\c(f,g)$ into a ``uniform $\Pi^0_1$ characteristic'' as defined
above.) \qedf\invdef

Some relations between these cardinal invariants are provable in ZFC: 
For example, if $g<  g' < f' < f$, then $\c(f',g')\le\c(f,g)$.   Also,
$\c(f^2, g^2) \le \c(f,g)$. 

 We will show 
  that if $(f,g)$ is sufficiently different from $(f', g')$, then the
  values of $\c(f,g)$ and $\c(f', g')$ are quite independent, and
moreover: if $\<(f_i, g_i):i<\omega_1>$ are pairwise sufficiently
different, then almost  any assignment of the form
$\c(f_i,g_i)=\kappa_i$ will be consistent. 

\bigskip

Similar results are possible for the  ``dual'' version of $\c(f,g)$:
$\c^d(f,g):=$ the smallest family of $g$-slaloms $\bar B$ such that
for  every $h$  bounded by $f$ there are infinitely many $k$ with
$h(k)\in B_k$, and for the ``tree'' version (a $g$-tree is a tree
where every node in level $k$ has $g(k)$ many  successors).  See
[\ShelahCDXLVIIIa]. 

			\bigskip

We thank Tomek Bartoszynski for pointing out the following known
results about the    cardinal characteristics $\c(f,g)$:  

  For example, lemma
   \constants{} follows from Theorem 3.17 in [\ComfortNegrepontis]:
Taking $\kappa = \alpha = \omega$, ${\beta} = n$, and letting ${\cal
S}\subseteq n^\omega$ be a family of $\omega$-large oscillation, then no 
family of $n{-}1$-slaloms of size $< 2^{\aleph_0}$  can cover $\cal
S$.   Indeed, whenever $F$ is a function on $\cal S$ such that for
each $s\in {\cal S}$, $F(s)$ is a $n{-}1$-slalom covering $s$, then $F$
has to be finite-to-one and in fact at most $n{-}1$-to-one.

Also, since  $\c(f,f{-}1)$ is the size of the smallest family of
   functions below $f$ which does not admit an ``infinitely equal''
   function, i.e., 
$$ \c(f,f{-}1)  = \min\{|G|: G \subseteq \prod_n f(n)\ \&\  
\forall h\in \prod _n f(n) \,\exists \, g\in G \,\forall^\infty n \,\,
   f(n)\not=g(n) \}$$
by [\Miller] we have that the minimal value of $\c(f,f{-}1)$ is the
smallest size of a set of reals which does not have strong measure zero. 

\smallskip
Also, note  that if $r$ is a random real over $V$ in $\prod_n f(n)$,
and     if $\sum_{n=1}^\infty 1/f(n)=\infty$, then $\prod_n (1 -
1/f(n)) =0$, so $r$ cannot be covered by any $f{-}1$-slalom from $V$.

Conversely, if  $\sum_{n=1}^\infty 1/f(n)<\infty$, then for any function
$h\in \prod_n f(n) \cut V$ there is a condition forcing that $h$ is
covered by the $f{-}1$-slalom  $(\{0,\ldots, f(k)-1\} - \{r(k)\} : k
\in \omega)$. 

Thus, if we add $\kappa$ many random reals with the measure algebra, a
easy density argument   shows that in the
resulting model we have 
$$ \c(f,f{-}1) = \cases{ \kappa = 2^{{\aleph_0}}& if 
         $\sum_{n=1}^\infty 1/f(n)=\infty$ \cr
\aleph_1 & otherwise (use any $\aleph_1$ many of the random reals)}$$

That already shows that we can have at least two distinct values of
$\c(f,g)$ and $\c(f',g')$.  

\bigskip

\noindent{\bf Contents of the paper:}  In section 1 we prove results
in ZFC of the form 
$$\hbox{``If $(f,g)$ is in relation \dots\  to $(f', g')$, 
	then $\c(f,g)\le \c(f',g')$''} 	$$

In section 2 we define a forcing notion $Q_{f,g}$ that
increases~$\c(f,g)$. (I.e., in $V^{Q_{f,g}}$, the $g$-slaloms from $V$ do not 
cover~$\prod_n f(n)$.)  Informally speaking, elements of $Q_{f,g}$ are
perfect trees  in which the size of the splitting is bounded by $f$,
sometimes $=1$, but often (i.e., on every branch), much bigger than~$g$. 

In section 3 we show that, assuming $\{(f_{\xi}, g_{\xi}):
{\xi}<\omega_1\}$ are sufficiently ``independent'', a countable
support product$\prod_{{\xi}<\omega_1} Q_{\xi}^{\kappa_{\xi}}$ of such
forcing notions will force $\forall {\xi}\, \c(f_{\xi}, g_{\xi}) =
\kappa_{\xi}$.

We use the symbol \SMILEY\ to denote the end of a proof, and we write
        \smiley\ when we leave a proof to the reader.

\bigskip\bigskip\bigskip

 \nextinput zfc {1. Results in ZFC}

\prop \ptwise Notation:  Operations and relations on functions are
understood to be pointwise, e.g., $f/g$, $g^\varepsilon$, $g<f$, etc.
$\[x]$ is the greatest integer~$\le x$.  $\lim f$ is $\lim_{k \to
\infty } f(k)$. 

We write $f\le^* g$ for $\exists n \, \forall  k\ge n \,\, f(k)\le
g(k)$.\qedf\ptwise

First we state some obvious facts: 

\prop \obvious Fact: 
\begindent
\ite 1 $f \le g$ iff~$\c(f,g)=1$. 
\ite 2 $f \le^* g$ iff $\c(f,g) $ finite. 
\ite 3 If $A:=\{k: g(k)< f(k) \}$ is infinite  then $\c(f\on A,	g\on
A) = \c(f,g)$. 
\ite 4  If $\pi $ is a permutation of $\omega$, then $\c(f\circ \pi,
	g\circ \pi) = \c(f,g)$.  \qed\obvious
\endent
(Strictly speaking, we define $\c(f,g) $ only for functions $f,g$
defined on all of $\omega$, so (3) should  be formally rephrased as
$\c(f\circ h, g \circ h) = \c(f,g)$, where $h$ is a 1-1 enumeration
of $A$)

\prop \obvRemark Convention: 
We will concentrate on the case where $\c(f,g)$ is infinite, so
	we will wlog assume that~$ g < f $.   By (4), we may also wlog
	assume that $g$ is nondecreasing.

In these cases we will have that $\c(f,g)$ is infinite, and moreover
an easy diagonal argument shows the following fact: 
\prop \uncount Fact: 
\begindent
\item{}
  $\c(f,g)$ is uncountable. \qed\uncount
\endent

Furthermore, we have the following properties: 
 
\prop \trivZFC Fact: 
\begindent
\ite 1 (Monononicity) If $f\le^* f'$, $g \ge^* g'$, then $\c(f,g) \le \c(f', g')$. 
\ite 2 (Multiplicativity)  $\c(f\cdot f',g\cdot g') \le \c(f,g)
		\cdot \c(f', g')$. 
\ite 3 (Transitivity)  $\c(f,h) \le \c(f, g) \cdot \c(g,h)$. 
\ite 4 (Invariance)  $\c(f,g) = \c(f^-, g^-)$ (where $f^-$ is the
	function defined by $f^-(n)=f(n+1)$. 
\ite 5 (Monotonicity II) If $A \subseteq \omega$ is infinite, then
       $\c(f\on A, g\on A) \le \c(f,g)$. \hfill\qed\trivZFC
\endent

\prop \trivZFCrem Remark: (2) implies in particular $\c(f^n, g^n) \le
\c(f,g)$.  See \squareex\ 
for an example of $\c(f^2, g^2) < \c(f,g)$. 

\bigskip

The following inequalities need a little more work.

\prop \moreZFC Lemma:
\begindent
\ite 1 $\c(f\cdot \left\[f /  g\right], f) = \c(f,g)$. 
\ite 2 $\c(f\cdot \left\[f /  g\right], g) =  \c(f,g)$. 
\ite 3 $\c(f\cdot \left\[f /  g\right]^m, g) = \c(f,g)$ for all
	$m\in \omega $.  
\endent

\proof   (2) follows from (1) using transitivity, and 
(3) follows from (2) by induction, so we only have to prove (1). 

Proof of (1):  By monotonicity we only have to show~$\le$.   So let $(N,{\in})$ be a
	reasonably closed model of a large fragment of ZFC (say,
	$(N,{\in}) < (H(\chi^+), {\in})$, where $\chi = 2^{\c}$) of size
	$\c(f,g)$ such $\prod_n f(n)$ is covered by the set of all
	$g$-slaloms from~$N$. 

Define $h$ by $h(k):= f(k)\cdot \left\[ f(k) / g(k)\right]$.  We can
	find  a family $\<B^i_k: i<f(k), k\in \omega>$ in $N$ such
	that for all $k$, $\{0, \ldots, h(k) -1 \} = \bigcup_{i<f(k) }
	B^i_k$, where $|B_k^i| \le {f(k) / g(k)}$.  We have to show
	that the set of $f$-slaloms from $N$ covers~$\prod_k h(k)$.

So let $x$ be a function satisfying $\forall
	k \, x(k)\in \bigcup _{i<f(k)} B_k^i$.  We can define a function
	$ y\in \prod_n f(n)$ such that for all $k$, $x(k)\in B_k^{y(k)}$.
	So there is some $g$-slalom $\bar C\in N$ such that for all $k$,
	$y(k) \in C_k$. 

Define $\bar A = \<A_k:k\in \omega > $ by $A_k := \bigcup_{i\in C_k} B_k^{i}$.
	Then $|A_k|\le |C_k | \cdot |B_k^i| \le g(k) \cdot {f(k) / 
	g(k)} = f(k)$, so $\bar A$ is an $f$-slalom in $N$, and  for all
	$k$,~$x(k)\in A_k$. 
\QED\moreZFC

\bigskip

\bigskip

\prop \techLemma Lemma: Assume~$f>g>0$. Assume that
	$\<w_i:i\in \omega>$ is a partition of $\omega$ into finite
	sets, and for each $i$ there are $\bar H^i = \<H^i_l:l\in
	w_i>$ satisfying (a)--(c). Then $\c(f', g') \le \c(f,g)$. 
\begindent
\ite a $\dom H^i_l = f'(i) = \{0, \ldots, f'(i)-1\}$
\ite b $\rng H^i_l \subseteq f(l)  = \{0, \ldots , f(l)-1\}$
\ite c Whenever $\<u_l:l\in w_i>$ satisfies 
\itemitem{} $u_l \subseteq f(l)$
\itemitem{} $|u_l| \le g(l)$
\item{} then $\{n< f'(i): \forall l\in w_i\, H^i_l(n)\in u_l\}$ has
	cardinality $\le g'(i)$
\endent

\proof  To any $g$-slalom $\bar B = \<B_l:l\in \omega>$ we can associate a
	$g'$-slalom $\bar B^* = \<B^*_i:i\in \omega>$ by letting $$
	B^*_i:= \{ n < f'(i): \forall l\in w_i\,\, H^i_l(n)\in w_l\}$$
Conversely, to any function $ x\in \prod_i f'(i)$
	we can define a function  $ x^*$ in $\prod_n f(n)$  by 
   $$ \hbox{if $l\in w_i$, then } x^*(l) = H^i_l(x(i))$$
It is easy to check that if $x^* $ is in $\bar B$  then  $x$ is
in~$\bar B^*$.  The result follows.  \QED\techLemma

\prop \productCor Corollary:  Assume $0=n_0 < n_1 < \cdots$, and let
$$\eqalign{ f'(i)&:= f(n_i) \cdot f(n_i+1) \cdots f(n_{i+1}-1)\cr 
	 g'(i) &:=
	g(n_i) \cdot g(n_i+1) \cdots g(n_{i+1} -1)}$$   

Then $\c(f', g') \le \c(f,g)$. 

\proof   Identify the set of numbers less than $f(n_i) \cdot f(n_i+1)
	\cdots f(f_{i+1}-1)$  with the cartesian product 
         $\prod_{n_i\le k< n_{i+1}}  f(k)$, and let 
         $$ H^i_l:\prod_{n_i\le k< n_{i+1}}  f(k) \ \to \ f(l)$$ 
	be the projection onto the $l$-coordinate.   We leave the
	verification of \techLemma(c) to the reader. \qed\productCor
							
\medskip

\prop \const Lemma: If $g$ is constant, $f(k) \ge 2^k$, then $\c(f,g)
   = \c$. 

\proof   Let   $\forall k \, g(k)=n$, $f(k ) = 2^k$.   Assume that
$\prod_l \pre l 2 $
	can be covered by $<\c$ many $g$-slaloms. 

For any $\eta\in \pre \omega 2 $, the sequence $\bar\eta:= \<\eta\on l:
   l\in \omega>$ is  in $\prod_l \pre l 2$.  But any $g$-slalom can
   contain only $n$ many such $\bar \eta$, i.e. for any $g$-slalom $\bar
	B = \<B_l:l\in \omega>$ we
   have $$ \left|\left\{ \eta\in \pre \omega 2 : \forall l \  \eta\on
	l \in    B_l\right\}\right| \le m $$
Since there are continuum many $\eta$ we need continuum many $g$-slaloms
	to cover $\prod_l f(l)$ (or equivalently, $\prod_l \pre l 2$).
\QED\const

\prop \constants Lemma:  If  $f$ and $g$ are constant with $f>g$, then
$\c(f,g) = \c$. 

\proof  Using monotonicity wlog we assume that $f(k) = n+1$, $g(k) = n$
	for all~$k$.  We will use 
   \techLemma.  Let $\omega = \bigcup_{i\in \omega} w_i$ be a
   partition of $\omega$ where $|w_i| = n^{2^i}$. 

Let $f'(i) = 2^i$, $g'(i)=n$,  and let $\<H^i_l:l\in w_i> $ enumerate
   all functions from $2^i$ to~$n$. 

We plan to show $\c(f,g) \ge \c(f',g')$ (so $\c(f,g)=\c$ by \const). 
   We want to apply \techLemma, so fix a sequence $\<u_l:l\in
   w_i>$,where  $u_l \subseteq f(l)$ and~$|u_l| \le g(l)$. 

To show that the hypotheses of \techLemma\ are satisfied, fix $i_0$  and
	let 
$$ A:= \{x< f'({i_0}): \forall l\in w_{i_0}\, H^{i_0}_l(x)\in u_l\} $$
and assume $A$ has 	cardinality~$> g'({i_0}) = n$.  So let $x_0,
   \ldots, x_n$ be distinct elements of~$A$.  Let $H:f'({i_0}) \to n+1$ be
   a function satisfying 
    $$ \forall j \le n\ H(x_j) = j$$
$H$ is one of the functions $\{H^{i_0}_l:l\in w_{i_0}\}$, say $H=H^{i_0}_{l_0}$.
	Let $j_0\notin 
	u_{l_0}$, then also  $$ x_{j_0} \notin \{x< f'({i_0}):  H^{i_0}_{l_0}(x)\in
	u_{l_0}\} \supseteq A,$$ contradicting~$x_{j_0}\in A$.  \QED\constants

\prop \diverge Corollary:  If $f>g$, and $ \liminf_{k\to\infty} g(k) <
   \infty$, then~$\c(f,g) = \c$. 

\proof  This follows from \constants, using monotonicity and
monotonicity II. \qed\diverge

\medskip

We can now extend \moreZFC\ as follows:

\prop \power Theorem:  
 If for some $\varepsilon>0$, $g^{1+\varepsilon} \le f$, then
	 for all $n$,~$\c(f^n, g) = \c(f, g)$.

\proof   First we consider a special case: Assume that~$g^2 \le f$.
Then we get  
$$\c(f,g) \le \c(f^2, g) \le \c(f^2, f) \cdot \c(f,g) \le \c(f^2, g^2)
	\cdot \c(f,g) = \c(f,g)$$
Now we use this result on $(f,g)$, then on $(f^2, g)$, etc, to get 
 $$ \c(f,g) = \c(f^2, g) = \c(f^4, g) = \c(f^8, g) = \cdots $$ 
and use monotonicity to get the general result under the
assumption~$g^2\le f$.  

Now we consider the general case $g^{1+\varepsilon} \le f$: 
\hfil\break
   If $g$ does not diverge to infinity, we have		already  (by
	\diverge)~$\c(f,g)= \c$.  Otherwise we can 	find some 
		${\delta} >0$ such that for almost all $k$, 
	$$ { f(k)  \over  g(k) } \ge g(k)^{\delta}+1, $$
so $$ \left\[f(k)  \over g(k) \right] \ge g(k)^{\delta}$$
Now choose $m$ such that~$m \cdot {\delta} >1$.
Then~$\left\[f(k)/g(k)\right]^m\ge g$. 	  
By \moreZFC, $ \c(f\cdot \left\[ f /  g\right]^m,g)  =   \c(f , g)$ 
and so by monotonicity also $\c(f\cdot g, g) = \c(f,g)$. 
Since $g^2 \le f\cdot g$, we can apply the result from the special
case above to get  $\c(f,g)= \c(f^n\cdot
		g^n, g) $ so in  particular, $\c(f^n, g) = \c(f,g)$. 
\QED\power

If $f$ is not much bigger than $g$, the assumption in \moreZFC\ and
\power\ may be false.  For these cases, we can prove the following:

\prop \addition Lemma: 
\begindent
\ite 1 $\c(2f-g, f) = \c(f,g)$. 
\ite 2 $\c(2f-g, g) =  \c(f,g)$. 
\ite 3 $\c(f+m(f-g), g) = \c(f,g)$ for all
	$m\in \omega $.  
\endent

\proof  The proof is similar to  the  proof of \moreZFC.  Again we only
have to show (1).  Let $(N,{\in})$ be a
	reasonably closed model of a large fragment of ZFC (say,
	$(N,{\in})\prec (H(\chi^+), {\in})$, where $\chi = 2^{\c}$) of size
	$\c(f,g)$ such $\prod_n f(n)$ is covered by the set of all
	$g$-slaloms from~$N$. 

Define $h$ by $h(k):= f(k) + f(k) - g(k)$.  We can
	find  a family $\<B^i_k: i<f(k), k\in \omega>$ in $N$ such
	that for all $k$, $\{0, \ldots, h(k) -1 \} = \bigcup_{i<f(k) }
	B^i_k$, where $|B_k^i| = 2$ for $l<f(k)-g(k)$, and $|B_k^i|=1$
        otherwise.   We have to show
	that the set of $f$-slaloms from $N$ covers~$\prod_k h(k)$.

So let $x$ be a function satisfying $\forall
	k \, x(k)\in \bigcup _{i<f(k)} B_k^i$.  We can define a function
	$ y\in \prod_n f(n)$ such that for all $k$, $x(k)\in B_k^{y(k)}$.
	So there is some $g$-slalom $\bar C\in N$ such that for all $k$,
	$y(k) \in C_k$. 

Define $\bar A = \<A_k:k\in \omega > $ by $A_k := \bigcup_{i\in C_k} B_k^{i}$.
	Thus $A_k$ is the union of $g(k)$ many sets, of which at most
        $f(k)-g(k)$ are pairs, and the others singletons.  Thus $|A_k|
        \le g(k) + (f(k)-g(k)) = f(k)$,  so $\bar A$ is an $f$-slalom in
        $N$, and  for all 	$k$, $x(k)\in A_k$. \QED\addition

\medskip

Similar to the proof of \power\ we now get: 

\prop \multiple Lemma: 
\begindent
\ite 1 If $2g \le f$, then for all $n$, $\c(nf, g) = \c(f,g)$. 
\ite 2 If for some $\varepsilon>0$, $(1+\varepsilon)g \le f$, then
	 for all $n$, $\c(nf, g) = \c(f, g)$. \hfill\qed\multiple
\endent

 \bigskip\bigskip\bigskip

 \nextinput definitions {2. The forcing notion $Q_{f,g}$}

\prop \nstarDef Definition:  We fix  sequences $\<n^-_k:k\in
\omega>$ and $\<n^+_k:k\in \omega> $ that increase very quickly
and satisfy $n^-_0 \ll n^+_0  \ \ll \ 
 n^-_1 \ll n^+_1 \ \ll \ \cdots$. 
 In particular, we demand
\begindent
\ite 1 For all $k$ $\displaystyle\prod_{j<k} n^-_j \le n^-_k $
\ite 2 $\displaystyle \lim_{k\to \infty}{ \log n^+_{k} \over
		\log n^-_k } = 0 $. 
\ite 3 $n^-_k\cdot n^+_k < n^-_{k+1} $. 
\endent\qedf\nstarDef

We will only consider functions $f$, $g$ satisfying $n^-_k \le g(k) <
f(k) \le n^+_k$.   This is partly justified by \productCor, and it
also helps to keep the formulation of the main theorem relatively
simple.

\prop \normdef Definition:  Let $X\not=\emptyset$ be finite, $c,d\in
\omega$.  A 
	$(c,d)$-complete norm on $\P(X)$ is a map 
$$\|\ |: {\bf P(X)}-\{\emptyset\}	\to	 \omega  $$
mapping any nonempty $a \subseteq X$ to a number $\|a|$ 
such that 
\begindent
\item{} whenever $a = a_1 \cup \cdots \cup a_c \subseteq X$, then for
	some $i_1, \ldots, i_d\in \{1,\ldots, c\}$, $\|a_{i_1}
\cup\cdots \cup  a_{i_d}|\ge \|a| -1$. ($|a|$ is the cardinality of
the set $a$)
\endent
      
A natural $(c,d)$-complete norm is given by $\|a| := \log_{c/d}|a|$.
$c$-complete means $(c,1)$-complete. \qedf\normdef

\bigskip\bigskip

\prop \thingydef Definition: We call $(f,g,h)$ {\bf \thingy}, if $f$, $g$,
$h$ are       functions in $\fct$,
 satisfying 
\begindent
\ite 1 For all $k$, $n^-_k \le  g(k) < f(k) \le  n^+_k$
\ite 2 For all $k$, $n_k^- \le h(k)$
\ite 3 $\lim_k \left.\log {f(k)\over g(k)}\right/ \log h(k) = \infty$. 
\endent

We call $(f,g)$ \thingy, if there is a function $h$ such that $(f,g,h)$
is \thingy  (or equivalently, if $(f,g, n^-)$ is \thingy,  where $n^-$ is
the function defined by $n^-(k)=n^-_k$). \qedf\thingydef

\bigskip\bigskip

\prop \fggRem Remark: For example, if $f$ and $g$ satisfy (1), then
$(f,g,g) $ is progressive iff $\log f / \log g \to \infty$. \qed\fggRem

\bigskip

In \singledef\ we will define a forcing notion $Q_{f,g,h} $ for any
\thingy~$(f,g,h)$.  First we recall the following notation:

\bigskip
\prop \Notation Notation: $\finSeq = \bigcup_n \pre n 2$ is the set of
finite sequences of 
natural numbers. For $s\in \finSeq$, $|s|$ is the length of $s$. 

A tree $p$ is a nonempty subset of $\finSeq$
with the properties
\begindent
\item{} $\forall  \eta\in p \,\forall k <|\eta|: \eta\on k\in p$
\item{} $\forall  \eta\in p : \sux_p(\eta) \not= \emptyset$, where
	$$\sux_p(\eta):= \{\nu \in p: \eta \subset \nu,
	|\eta|+1=|\nu|\}.$$
\endent
A branch $b$ of $p$ is a maximal  linearly $\subseteq $-ordered subset
of~$p$.  Every branch $b$ defines a function $\bar b:\omega\to \omega$
by $\bar b = \bigcup b$.  We usually identify $b$ and $\bar b$, so we
write $b\on k$ (instead of $(\bigcup b)\on k$) for the $k$th 
element of~$b$.   

The set of all branches of $p$ is written as~$[p]$. 

For $\eta\in p$, we let 
$$ \trim p^\eta := \{\nu\in p :  \nu \subseteq \eta \hbox{ or } \eta
			\subseteq \nu\}$$

We let 
$$\eqalign{ \split(p)&:=\{\eta\in p: |\sux_p(\eta)| >1 \}
\qquad \hbox{(the splitting nodes of $p$)}
\cr
\split_n(p) &:= \{ \eta\in \split(p): |\{\vu\subset \eta: \vu\in
\split(p)\}| = n \} \quad \hbox{(the $n$-th splitting level)}
\cr}$$
and we define the stem  of $p$ to be the unique element of
$\split_0(p)$.  
\qedf\Notation
\bigskip

\prop\singledef Definition:  Assume $f,g, h$ 
	are as in \thingydef. Then we define for all $k$, and for all
sets $x$ 
 $$ \| x |_k \ := \ \left\[ \log(|x|/g(k)) \over \log h(k)\right]$$
and we define 
 the forcing notion
	$Q_{f,g}$ (or more accurately, $Q_{f,g, h}$) to be the
	set of all $p$ satisfying 

\begindent
\ite 1 $ p$ is a perfect tree. 
\ite 2 $\forall \eta\in p \, \forall i\in\dom(\eta)\,\, \eta(i) <
	f(i)$. 
\ite 3 $\forall \eta\in \split_n(p) \|\sux_p(\nu)|_{|\nu|} \ge  n$. 
\endent

We let $p \weaker q $ (``$q$ extends $p$'')  iff $q \subseteq p$. \qedf\singledef

\prop \AxiomA Remark: If we define 
$$ p \sqsubseteq_k q  \hbox{ iff } p \weaker  q \hbox{ and } \split_k(p)
\subseteq q$$
then $Q_{f,g,h}$ satisfies axiom A, and is in fact strongly
$\fct$-bounding, i.e., for name of an ordinal, $\nalpha$, for any $p$
and for any $n$ there is a finite set $A$ and a condition $q
\sqsupseteq_n p$, $q \forces \nalpha \in A$.   However, it will be
more convenient to use  the relation $\le_n$ that is based on {\it
levels} rather than {\it splitting~levels}. \qedf\AxiomA

\prop \leqsubn Definition:  For $p$, $q\in Q$, $n\in \omega$  we define 
$$ p \weaker_n q \ \hbox{ iff } \ p \weaker q \hbox{ and } p\cut
\pre{\le n}\omega \subseteq q\eqno \qedf\leqsubn$$

\prop \normnotation Notation:  We will usually write $\|\eta|_p$
	instead of $\|\sux_p(\eta)|_{|\eta|}$. \qedf\normnotation

\prop \cfCCCXXVI Remark:  This forcing is similar to the forcing in
	[\ShelahCCCXXVI], but note the following important difference:
	Whereas in [\ShelahCCCXXVI] all nodes above the stem have to
	be splitting points, we allow many nodes to have only one
	successor, as long as there ``many''   nodes with
	high norm. \qedf \cfCCCXXVI

\prop \normrem Remark: 

\begindent
\ite 1 The norm $\|\cdot|_k$ is $h(k)$-complete (hence
	also $n^-_k$-complete). 
\ite 2 If $c/d\le h(k)$, then the norm is $(c,d)$-complete. 
\ite 3 If $\|a|_k > 0$, then $|a| > g(k)$. 
\ite 4 $\|f(k)|_{k} \to \infty$ (so $Q_{f,g,h}$ is nonempty). \qed \normrem
\endent

We will see  in the next section that this forcing (and any countable
support product of such forcings) is proper and
$\fct$-bounding.   For the moment, we only show why this forcing is
useful in connection with $\c(f,g)$: 

\prop \olddontcover Fact: Any generic filter  $G \subseteq Q_{f,g}$
defines a ``generic branch'' $$ r := \bigcup_{p\in G} \stem(p) $$
that avoids all $g$-slaloms from~$V$. 

\proof  Let $\bar B = \<B_k:k\in \omega > $ be a $g$-slalom in $V$, and
let $p\in Q_{f,g}$ be a condition. Let $\eta\in p$ be a node
satisfying $\|\eta|_p > 0$. Let $k:=|\eta|$. 
 Then $|\sux_p(\eta) | > g(k)$ by  \normrem(3), so there is $i\notin B_k$,
$\eta\extend i \in p$. So $\trim p ^ {\eta\extend i } \forces
r(k)=i\notin B_k$. \QED \olddontcover

\bigskip\bigskip\bigskip

\nextinput construction {3. The construction}

In this section we will prove the following theorem:

\prop \mainTHEO Theorem (CH): Assume that $(f_{\xi},g_{\xi} : 
{\xi}< \omega_1)$ is a sequence of \thingy\ functions, witnessed by
functions $h_{\xi}$ (see \thingydef).

 Let $(\kappa_{\xi}:{\xi}<\omega_1)$ be a
  sequence of      cardinals satisfying
  $\kappa_{\xi}^\omega=\kappa_{\xi}$ such that  	 
     whenever $ \kappa_{\xi} < \kappa_\zeta $, then 
$$ \lim_{k\to \infty} \min\left( {f_{\zeta} (k)\over g_{\xi}(k)}, 
            \left.{f_{\xi}(k)\over g_{\xi}(k)}\right/ h_{\zeta}(k)
\right) \ = \ 0 $$
(or informally: either $f_{\zeta} \ll g_{\xi}$, or $f_{\xi}/g_{\xi}
\ll h_{\zeta}$, or a combination of these two condition holds)

Then there is a proper forcing notion $P$ not collapsing cardinals
     nor changing cofinalities such that $$\forces_P \forall {\xi}:
     \c(f_{\xi}, g_{\xi}) = \kappa_{\xi}$$

\medskip

For the proof we use a countable support product of the forcing
notions $Q_{f_{\xi}, g_{\xi}, h_{\xi}}$ described in the previous
section.   

\prop \onlytwo Remark: 
The theorem is of course also true (with the same proof) if we have
countably or finitely 
many functions to deal with.

If we are only interested in 2 cardinal invariants $\c(f',g')$, $\c(f,g)$,
then we can phrase the theorem without the auxiliary functions $h$ as
follows:  If $(f,g)$ and $(f', g')$ are \thingy, and satisfy 
$$  \min \left( {f' \over g}, { \log(f/g) \over \log (f'/g')}\right) \
\to \ 0$$
 then $\c(f,g) < \c(f', g') $ is consistent. 

In particular, this shows that our result is quite sharp:  For
example, if for some function $d$ we have $\lim d = \infty$,  $f' = f^d$,
$g'=g^d$ (and $(f,g)$, $(f', g')$ are progressive with the same
$n_k^-$, $n_k^+$), then $\c(f,g) < \c(f', g')$ is consistent.  On the
other hand, $\c(f^n, g^n) \le \c(f,g)$ for every fixed $n$.

\proof  Choose $h'$ such that $\log h'  \approx 2 \log (f/g)$
whenever ${f' \over g} \ge  { \log(f/g) \over \log (f'/g')}$. 
$(f',g',h')$ is progressive, and the assumptions of the theorem are
	satisfied.  (Recall that  $(f,g)$  is progressive, hence $\log
        f/g\gg \log n^-$, so $h'$ will satisfy $h'(k) \ge n^-_k$).
        \QED \onlytwo

A similar simplified formulation of \mainTHEO{} is possible when we
deal with only countably   many functions.

\bigskip

\prop\Blassanswer Example: There is a family $\<(f_{\xi}, g_{\xi},
g_{\xi}: {\xi} < \c> $ of continuum many progressive functions such
that for any $\zeta \not= {\xi}$, $ \min \left({f_{\xi}\over
g_{\zeta}}, {f_\zeta \over g_{\xi}}\right) \to 0$.   [In particular,
under CH we may choose any family $(\kappa_{\xi}: {\xi} < \omega_1)$
of cardinals satisfying $\kappa_{\xi}^\omega = \kappa_{\xi}$ and get
an extension where $\c(f_{\xi}, g_{\xi})= \kappa_{\xi}$.]

\proof  Let $\ell_k:= \left\[ {1\over 2} \sqrt{\log {\log n_k^+\over
\log n_k^-}}\right]$.  (Here, ``$\log$'' can be the logarithm to any
(fixed) base, say 2.)    Then $\lim_{k\to \infty}\ell_k = \infty$, and by
invariance (\trivZFC(4)) we may assume $\ell_k \ge 1 $ for all $k$. 

Let $T \subseteq 2^{{<}\omega}$ be a perfect tree such that for all
$k$ we have $|T \cut 2^k| = \ell_k$, say, $T\cut 2^k = \{s_1(k),
\ldots, s_{\ell_k}(k)\}$.  

For any $x \in [T]$
 (i.e., $x\in 2^\omega$, $\forall k \, x\on k \in T$) we now define
functions $f_x$, $g_x$, $h_x$ by:  
\begindent
\item {}  If $x\on k = s_i(k)$, then 
$$\eqalign{ f_x(k) &= \left(n_k^-\right)^{\ell_k^{2i}} \cr
      h_x(k)= g_x(k) &= \left(n_k^-\right)^{\ell_k^{2i-1}} \cr}$$
\endent
We leave the verification that $(f_x, g_x, h_x)$ is indeed progressive
to the reader. [Recall 2.4, and also note that  $\log \log f_x(k) \le
2 \ell_k \log \ell_k + \log \log n^-_k < \log\log n^+_k$. 
Finally, note  that if $x \not= y$, then for almost all $k$ we
have $\displaystyle \min \left( {f_x(k)\over g_y(k)}, {f_y(k) \over
h_x(k)} \right) \ll  {1 \over n_k^-}$.]  \QED\Blassanswer

\prop\squareex Example:  It is consistent to have $\c(f^2, g^2) <
\c(f,g)$ (for certain $f$, $g$).  

\proof  Let $\ell_k:= \displaystyle\left\[ {1 \over 6} \log {n_k^+\over n_k^-}
\right]$.  Assume $\ell_k > 0$ for all $k$.   Then, letting 
    $$\eqalign {f(k) &:= \left(n_k^-\right)^{3\ell_k}\cr
                g(k) &:= \left(n_k^-\right)^{2\ell_k}\cr
                h(k) &:= n_k^-\cr}$$
We have that $(f,g, h) $ and $(f^2, g^2, h)$ are progressive, and
$\lim {f \over g^2} = 0$, so we can apply  the theorem. \QED\squareex

\prop \iterdef Definition:  

Let $\kappa $ be 
  a disjoint union $\kappa = \bigcup_{{\xi}<\omega_1} A_{\xi}$,
where $|A_{\xi}|= \kappa_{\xi} $. 
  
For $ \alpha < \kappa $, let $Q_\alpha $ be the forcing $Q_{f_{\xi},
  g_{\xi}, h_{\xi}} $, if $\alpha \in A_{\xi} $, and let $P =
  \prod_{\alpha<\kappa}Q_\alpha$ be the {\bf countable support
	product}  of   the forcing notions $Q_\alpha$, i.e., elements
	of $P$ are countable functions $p$ with $\dom(p) \subseteq
	\kappa$, and $\forall \alpha\in \dom(p) \,\, p(\alpha)\in
	Q_\alpha$. 

For $A \subseteq \kappa$, we write $P\on A:= \{p\on A: p\in P\}$.
  Clearly $P \on A \lesscirc P$ for any~$A$.  In particular,~$Q_\alpha
  \lesscirc P $.

We write $\name r_\alpha$ for the $Q_\alpha$-name (or $P$-name)  for the
  generic branch introduced by a generic filter on~$Q_\alpha$. 

We say that $q $ {\bf strictly extends} $p$, if $q \stronger p $ and
  $\dom (q) = \dom(p)$. \qedf\iterdef

\prop\facts Facts: Assume CH.  Then 
\begindent
\ite 1 each $Q_\alpha$ is proper and $\fct$-bounding.
\ite 2 $P$ is proper and $\fct$-bounding. 
\ite 3  $P$ satisfies the $\aleph_2$-cc. 
\ite 4  Neither cardinals nor cofinalities  are changed by forcing
		with~$P$.  
\endent

Proof of (1), (2):  See below (\propercoro, \properremark)

Proof of (3):  A straightforward ${\Delta} $-system argument, using
CH. 

(4) follows from  (2) and (3). 
\QED\facts

We plan to show that ${} \forces_P \c_{\xi} = \kappa_\xi$ for all
  ${\xi} < \omega_1$.

\prop \leveldef Definition:  If $p\in P$, $k\in \omega$, we
let the level $k$ of $p$ be $$\level_k(p) :=
\theset \bareta : $\dom({\bareta}) = \dom(p)$, \hskip 3cm \\
   $\forall \alpha \in \dom({\bareta}): |{\bareta}(\alpha)|=k$,
${\bareta}(\alpha)\in p(\alpha)$
  \endset
$$
We define the set of active ordinals at level $k$ as 
$$\act_k(p) := \{\alpha\in \dom(p): |\stem(p(\alpha))| \le k \}$$

\prop \levelrem Remark: Sometimes we identify the set $\level_k(p) $
with the set 
$$\eqalign{&\{{\bareta} : \dom({\bareta}) = \act_k(p),  \forall \alpha \in
\dom({\bareta}): |{\bareta}(\alpha)|=k\}\cr
&\quad = \{ {\bareta} \on \act_k(p): {\bareta} \in \level_k(p)\}}$$

\prop \splitlevel Definition:  We say that  the $k$th level
is a splitting level of $p$ (or ``$k$ is a splitting level of $p$'') iff 
$$ \exists \alpha \in \dom (p) \,\exists \eta\in \split(p(\alpha)):
	|\eta| = k$$

\prop\extendDf Definition:  If ${\bareta} \in\level_k(p)$, ${\bareta} '
  \in \level_{k'}(p)$, $k<k'$,  then we say that ${\bareta} ' $
  extends ${\bareta}$ iff for all $\alpha \in \dom({\bareta})$,
  ${\bareta}'(\alpha) $ extends (i.e., $\supseteq$)~${\bareta}(\alpha)$.

\prop \fkDef Definition:  For $p,q\in P$, $k\in \omega$, 
we let 
$$ p \weaker_{k} q \ \hbox{ iff } \ p \weaker q \hbox{ and } 
	\forall \alpha\in  \dom(p):  p(\alpha)\weaker_k q(\alpha)
\hbox{ and } \act_k(p) = \act_k(q)$$
That is, we allow $\dom(q)$ to be bigger  than $\dom(p)$, but for all
new $\alpha \in \dom(q)-\dom(p)$ we require that $|\stem(q(\alpha))| > k$.

\prop\densedef Definition:  Let~$A \subseteq P$.  A set $D \subseteq P$ is 
\begindent
\advance\parindent by 3cm\advance\rightskip by -2cm\relax
\item{dense in $A$,} if $\forall p\in A\, \exists q \in D : p \weaker q$
\item{strictly dense in $A$,} if  $\forall p\in A\, \exists q \in D :
	p \weaker q$ 
	and $\dom (p) = \dom(q)$
\item{open in $A$,} if $\forall p \in D\,\forall q\in A$: ($p \weaker q$ implies
	$q\in D$)
\item{almost open in $A$,} if  $\forall p \in D\,\forall q\in A$: ($p
\weaker q$ 	and $\dom(p)=\dom(q) $ implies 	$q\in D$)
\endent
      These definitions can also be relativized to conditions above a
	given  condition~$p_0$.    If we omit $A$ we mean~$A=P$.

\prop\bartrim Definition:  If ${\bareta} \in \level_k(p)$, we let
	$q=\trim p ^{\bareta} $ be the condition defined by 
	$ \dom(q)=\dom(p)$, and 
 $$ \forall \alpha \in \dom(q) \,\, q(\alpha) = \trim
	p(\alpha)^{\bareta(\alpha)} $$

\prop \decidedef Definition:  If $p \forces \name x \in V$, and
	${\bareta} \in \level_k(p)$, we say that ${\bareta}$ decides
	$\nx$ (or more accurately, $\trim p^{\bareta} $ decides $\nx$)
	if for some $y\in V$, $\trim p^ {\bareta}  \forces \nx =
	\check y$.

First we simplify the form of our conditions such that all levels are
finite.

\prop \activefinite Fact: The set of all conditions $p$ satisfying 
\begindent
\frame{finite} $\forall k\, |\act_k(p)|<\omega $, and moreover: 
\frame{onepair} For any splitting level $k$ there is exactly one pair
		$(\eta,\alpha) $ such that
		$|\sux_{p(\alpha)}(\eta)|>1$.  
\endent
is dense in $P$. \qed \activefinite

\prop\levelsize Fact: If $p$ is in the dense set given by \:{finite}
	and \:{onepair}, then 
        the size of level $k$ is $\le n_{k-1} ^-\cdot n_{k-1}^+ < n_k^-$. 

\proof  By induction. \qed\levelsize

{\bf From now on we will only work in the dense set of conditions 
satisfying \:{finite} and \:{onepair}.}

\prop \kldef Notation:  For $p$ satisfying \:{finite}--\:{onepair},
	 we let  	$k_l= k_l(p) $ be the $l$th splitting
	level. 
Let $\eta_l = \eta_l(p) $ and $\alpha_l = \alpha_l(p) $ be such that
	$|\eta_l(p)| = k_l(p)$, $\eta_l(p)\in \split(p(\alpha_l))$. 
    We let ${\zeta}_l = {\zeta}_l(p) $ be such that $\alpha_l \in
    A_{{\zeta}_l}$.

We write $\|p|_{k_l}$ for $\|\eta_l|_{p(\alpha_l)}$, i.e., for
$\|\sux_{p(\alpha_l)}(\eta_l)|_{\zeta_l, k_l}$. (See figure 2)

\midinsert
$$\vcenter{\epsffile[35 97 270 298]{448-pic1.ps}}$$
\centerline{Figure 2: A condition satisfying (I) and (II)}
\endinsert

\prop \pruneDef Definition:  
If $p$ is a condition, $l\in \omega$, $\alpha^*:=\alpha_l(p)$,
	$\eta^*:=\eta_l(p)$, $\nu^*\in \sux_{p(\alpha^*)}(\eta^*)$, we can
	define a stronger condition $q$ by letting
	$q(\alpha)=p(\alpha)$ for all $\alpha\not= \alpha^*$, and 
 $$ q({\alpha^*}) := \{\eta\in p(\alpha^*): \hbox{If $\eta^* \subset
	\eta$, then $\nu^* \subseteq \eta$}\}$$
In this case, we say that $q$ was obtained from $p$ by ``pruning the
	splitting node $\eta^*$.''

To simplify the notation in  the fusion arguments below, we will use the
following game:

\prop \gamedef Definition:  For any condition $p\in P$, $G(P,p)$ is
the following two person game with perfect information: 

{\advance\leftskip by 0.3cm
\advance\rightskip by 0.3cm\relax

There  are two players, the \spendthrift/ and the \accountant/.  A play in
$G(P,p)$ last $\omega$ many moves (starting with move number 1) The
\accountant/ moves first. We let $p_0:=p$,~$i_0:=0$. 

In the $n$-th move, the \accountant/ plays a pair $(\eta^n, \alpha^n)$ with
$\eta^n\in p_{n-1}(\alpha^n)$, $|\eta^n|=i_{n-1}$, and a number~$b_n$. 

Player \spendthrift/ responds by playing a condition $p_{n}$ and a
finite sequence $\nu^n$ (letting $i_n:=|\nu^n|+1$) satisfying the
	following:  (See Figure 3)
\begindent
\ite 1 $p_n \stronger_{i_{{n-1}}} p_{n-1} $.
\ite 2 $\nu^n \in p_n(\alpha^n)$
\ite 3 $\|\vu^n|_{p_n(\alpha^n)} > b_n$.
\ite 4 $\nu ^n \supset \eta^n$. 
\ite 5 For all $\alpha\in \dom (p_n) - \dom (p_{n-1} )$,
$|\stem(p_n(\alpha))| > |\nu^n|$. 
\ite 6 $|\level_{|\vu^n|}(p_n)| = |\level_{|\eta^n|}(p_n)| 
	= |\level_{|\eta^n|}(p_{n-1} )| $
\endent
      
(Remember that all conditions $p_n$  have to be in the dense
	set given by \:{finite} and \:{onepair}) 
Player \accountant/ wins iff after $\omega$ many moves there is a
	condition $q$ such that for all $n$, $p_n \weaker q$, or
equivalently, if the function $q$ with domain $\bigcup_n \dom (p_n)$,
defined by 
$$ q(\alpha) = \bigcup _{\alpha \in \dom (p_n)} p_n(\alpha)$$
is a condition. Note that we have $\eta_l(q) = \nu^l$, $\alpha_l(q) =
	\alpha^l$, since the only splitting points are the ones
	chosen by \spendthrift/. 

}
\midinsert
$$\vcenter{\epsffile[9 79 297 317]{448-pic2.ps}}$$
\centerline{Figure 3: stage $n$}
\endinsert
\prop \strategy Fact:  Player \accountant/ has a winning strategy in~$G(P,p)$.

\proof  We leave the  proof to the reader, after pointing out
	that a finitary bookkeeping will ensure that the limit of the
	conditions $p_n$ is in fact 	a condition. 

In particular, this shows that \spendthrift/ has no winning strategy.
	Below we will define various strategies for the \spendthrift/,
	and use only the fact that there is  a  play in which the
	\accountant/ wins. \qed \strategy

\bigskip 

The game gives us the following lemma: 

\prop \thinning Lemma: Assume that $p$ is a condition satisfying
    \:{finite}--\:{onepair}. 
For each $l$ let $\emptyset \not= F_{\eta_l} \subseteq
  \sux_{p(\alpha_l)}(\eta_l)$ be a set 
of norm $\| F_{\eta_l}|_{k_l} \ge \| \sux_{p(\alpha_l)}(\eta_l)|/2$. 

Then there is a condition $q \stronger p $, $\dom (q) = \dom(p)$ such
	that for all $l$:  
$$ (*)\qquad\qquad\hbox{ If $\eta_l(p)\in q(\alpha_l(p))$, then 
 $\sux_{q(\alpha_l(p))}(\eta_l(p)) \subseteq F_{\eta_l}$}$$

\proof   The condition $q$ can be constructed by playing the game.
In the $n$-th move, \spendthrift/ first finds a  $\eta^n \supset \nu^n$ 
satisfying  $\eta^n(i) \in F_{\eta_i}$ whenever this is applicable, 
and 
$\| \sux_{p_{n-1}}(\eta^n)| > 2 b_n$. 
Then \spendthrift/ obtains $p_n$ by pruning  (see \pruneDef) all splitting nodes of
	$p_{n-1}$ whose 
	height  is between $|\eta^n|$ and $|\nu^n|$ and further
	thinning out the successors of $\eta^n$ to
	satisfy~$\sux_{p_n}(\eta^n) = F_{\eta^n}$.  (Note that
	$F_{\eta^n} \subseteq \sux_{p_{n-1} }(\eta^n) = \sux_{p_0
	}(\eta^n)$.) 

In the resulting condition $q$ the only splitting nodes will be the
	nodes~$\eta^n$, so $(*)$ will be satisfied. \QED \thinning

(Note that in general   $\eta_l(q) \not= \eta_l(p)$, and indeed
  $k_l(q) \not = k_l(p)$, since many splitting levels of $p$ are not
  splitting levels in $q$ anymore.)

\prop \densefour Lemma:  Assume $\ntau$ is a $P$-name of a function
from $\omega$ to~$\omega$, or even from $\omega$ into ordinals.  Then
	the set of         conditions satisfying   
	\:{finite}--\:{kplusone}  is 	dense 
	and almost open.
\begindent
\frame{kplusone} Whenever $k$ is a splitting level, then every
	${\bareta} $ in 
	level ${k+1} $ decides~$\ntau \on k$.
\endent
      
      Proof of \:{kplusone}:  We will use the game from \gamedef. We
	will define a strategy for the \spendthrift/ ensuring that the
	condition $q$ the \accountant/ produces at the end will satisfy
	\:{kplusone}.    

 In the $n$-th move, \spendthrift/ finds a condition $r_n
	\stronger_{i_{n-1}} p_{n-1}$ such that for every ${\bareta}\in
	\level_{i_{n-1}}(r_n)$ the condition $\trim {(p_n)}^ {\bareta}
	$ decides~$\ntau \on i_{n-1} + 10$.  Then \spendthrift/ finds
	$\eta^n\in r_n(\alpha^n)$ satisfying the rules and obtains
	$p_n$ with $\eta^n\in p_n(\alpha^n)$ from $r_n$ by pruning all
	splitting levels between $i_{n-1}$ and~$|\eta_n|$.
\QED \densefour

Since all levels of $q$ are finite, it is thus possible to find a
	finite sequence $\bar B = \<B_k:k\in \omega>$  in the ground
	model that will cover 	$\ntau$.  (I.e.\ $q \forces \ntau(k)
	\in B_k$).      The rest of this section
	will be devoted to finding ``small'' such sets $B_k$.

\prop\propercoro Corollary:   $P$ is $\fct$-bounding and does not
collapse $\omega_1$. \qed \propercoro

\prop \properremark Remark:  
Although it does not literally follow from 	\densefour, the reader
	will have no difficulty in showing that 	$P$ is
	actually $\alpha$-proper for any $\alpha< \omega_1$. \smiley\
        Indeed,
	using the partial orders $\sqsubseteq_n$ from \AxiomA, it is
	possible to carry out straightforward fusion arguments,
	without using the game \gamedef\ at all.
	However, the orderings  $\le_n$ are more
	easy to handle, since in induction steps we
	only have to take care of a single
	$\eta^n$, instead of a front.

\prop \smallsupport Fact: $\forces_P \forall \tau\in\fct\, \exists B
	\subseteq  \kappa$, $B$ countable, $B\in V$, and~$\tau \in
	V[G\on B]$.  

\proof   Let $p$ be any condition and let $\ntau$ be a name for a real.
There is a stronger condition $q$ satisfying \:{finite}, \:{onepair}
and \:{kplusone}.   Let~$B:=\dom(q)$.
Clearly~$q \forces \ntau \in V[G\on B]$. \QED\smallsupport

\prop\supportcor Corollary:  If ${\lambda} = |A|^\omega$,
	then~$\forces_{P\on A} {2^{\aleph_0}} \le {\lambda}$.  

\proof   For each countable subset $B \subseteq A$, $\forces _{P\on B}
	CH$.   Since every real in $V[G]$ is in 
	some such $V[G\on B]$, the result follows.  \QED\supportcor

\prop \plusnotation Fact and Notation: If $p$ satisfies \:{onepair},
then 
\begindent
\ite 1 If ${\bareta} (\alpha_l) =\eta_l$, and $\nu\in
	\sux_{p(\alpha_l)}(\eta_l)$,  then the requirement 
        $$ {\bareta} ^{+\nu}(\alpha_l)=\nu$$
        uniquely defines an extension $ {\bareta} ^{+\nu}$
	of ${\bareta}$ in  $ \level_{k_l+1}(p)$. 
\ite 2 If ${\bareta} (\alpha_l) \not=\eta_l$,  ${\bareta}$ has a
	unique extension ${\bareta}^{+} \in \level_{k_l+1}(p)$.  To
	simplify the notation in \leproof{} below, we  also define for
	this 	case,  for any $\vu \in
	\sux_{p(\alpha_l)}(\eta_l)$, ${\bareta}^{+\nu} :=
	{\bareta}^+$. 
\endent

\prop\smalll Fact: The set of conditions satisfying \:{smalllevel} is
	strictly dense  (but not almost open) in the set of conditions satisfying
        \:{finite}--\:{onepair}.   

\begindent
\frame{smalllevel} For all $l$:
 $$|\level_{k_{l}}(p)| < \min\left({\|p|_{k_l} \over 2}, n_{k_l}^-\right) $$
\endent
For the proof, note that 
$|\level_{k_{l}}(p)| = |\level_{k_{l-1}+1}(p)| $. \qed\smalll

\bigskip

\prop \densetwo Lemma:  Assume $\ntau$ is a $P$-name of a function
		$\in \fct$,
and $\forces _P \forall k \, \ntau(k)<n^+_{k}$.  Then the set of
conditions satisfying \:{tauonk} is strictly dense and almost open  in
	the set give by 	\:{finite}, \:{onepair}, \:{kplusone}.
	where
\begindent
\frame{tauonk}  Whenever $k$ is a splitting level, then every
	${\bareta} $ in 
	level $k$ decides~$\ntau\on k$. 
\endent

\proof  Fix $p$ satisfying 
	\:{finite}, \:{onepair}, \:{kplusone}, \:{smalllevel}.

  Let $k_l:= k_l(p)$, etc. 
  Let~$m_l:= |\level_{k_l}|$. 
%
%
%
\proof We will use \thinning.  For each $l\in \omega$, $F_{\eta_l} \subseteq
	\sux_{p(\alpha_l)}( \eta_l) $ will be defined as follows:  Let
	$m_l :=| \level_{k_l} (p) |$, and let ${\bareta}^0, \ldots,
	{\bareta}^{m-1}$ enumerate $\level_k(p)$.  Find a sequence 
 $$ 	\sux_{p(\alpha_l)}( \eta_l)  = F^0 \supseteq F^1 \supseteq
	\cdots \supseteq F^m\qquad \forall i \|F^{i+1}|_k \ge \|F^i|_k-1
	$$ 
such that for all $i$ there exists $x^i$ such that for all $\vu \in
	F^{i+1} $ we have $\trim p^{({\bareta}^i)^{+\nu}} \forces \ntau
	\on k = x$.   It is possible to find such $F^{i+1}$ since
	$\|\cdot|_k$ is $n_k^-$-complete, and there are only $n_0^+
	\cdot n_1^+ \cdots n_{k-1}^+ < n_k^-$ many possible values of
	$\ntau \on k$.

Finally, let $F_{\eta_l}:= F^m$.   Applying
	\thinning\ will yield the desired result.  \QED\densetwo

\prop \densesix Remark:  Note that \:{tauonk} in particular implies 
\begindent
\frame{tau(k)} Whenever $k$ is not a splitting level, then every
  ${\bareta} $ 
	in level $k$ decides~$\ntau(k)$. 
\endent
      
\bigskip

\prop \geproof Proof that $\forces_P \c(f_{\xi} , g_{\xi})\ge
	\kappa_{\xi}$: (This proof is essentially the same as
	\olddontcover.) 

Recall that  $\name r_\alpha$ is  the generic real added by the forcing~$Q_\alpha$.  
Working in $V[G]$, let $\B$ be a family of less than $\kappa_{\xi}$
many $g_{\xi}$-slaloms.  We will show that they cannot cover $\prod
	f_{\xi}$,  by finding an $\alpha$ such that $\name r_\alpha$
	is forced not to be covered. 

There exists a set
$A\in V$ of size $<\kappa_{\xi}$ such that $\B \subseteq V[G\on A]$.
	Since $|A|< \kappa_{\xi}$ there is $\alpha\in A_{\xi} - A$. 

Assume that  $\bar B$ is  a $g_{\xi}$-slalom in $V[G\on A]$ covering~$r_\alpha$.
So in $V$ there are  a $P\on A$-name $\name{\bar B}$ and a condition
$p$ such that 
$$ \forces_{P\on A} \name {\bar B} \hbox{ is a $g$-slalom}$$
and 
$$ p \forces _P \name {\bar B} \hbox{ covers $r_\alpha $}$$

We can find a node $\eta$ in $p(\alpha)$ with $\sux_{p(\alpha)}(\eta)$
having more than $g(|\eta|)$ elements. Increase $p\on A$ to decide
$\name B_{|\eta|}$, then increase $p(\alpha)$ to make $r_\alpha$ avoid
	this set.  
\QED \geproof

\prop \almostall Fact: 
Fix ${\xi^*}$.  Then the  set of conditions $p$ satisfying 
\begindent
\frame{almostall} For all $l$:  If 
$\kappa_{\xi^*}<\kappa_{\zeta_l(p)}$, then 
$$
 \min\left( {f_{{\zeta_l(p)}} ({k_l})\over g_{\xi^*}({k_l})}, 
            \left.{f_{\xi^*}({k_l})\over g_{\xi^*}({k_l})}\right/
	h_{{\zeta_l(p)}}({k_l}) \right) 
          \ < \ 
	{1 \over |\level_{k_l}(p)|}$$
\endent
is dense almost open. 

\proof   Write $F_{{\zeta} }$ for the function 
$\displaystyle \min\left( {f_{\zeta} \over g_{\xi^*}}, 
            \left.{f_{\xi^*}\over g_{\xi^*}}\right/ h_{\zeta} \right)$.
Recall that  if $\kappa_{\zeta} < \kappa_{{\xi}^*} $, then
	$F_{\zeta}$ tends to $0$. 

Fix a condition $p $, 
	We will use the game $G(P,p)$.  \spendthrift/ will use the
	following strategy:  Whenever $\alpha_n \in A_{{\zeta}}$ and
	$\kappa_{\zeta}<\kappa_{{\xi}^*}$ , then
	\spendthrift/ first find $m_0$ such that for all $m \ge m_0$ we
	have $F_{{\zeta} }(m) < 1/|\level_{h_{n-1}}(p_{n-1})|$.   Now find
	$\nu^n\supseteq \eta^n$ of length $> m_0$ with a large enough
	norm, and play any condition $p_n$ obeying the rules of the
	game.   In particular, we must have $\left|\level_{|\nu^n|}(p^n)\right| =
	\left|\level_{|\eta^n|}(p^n)\right|$. 

Clearly the condition resulting from the game satisfies the
requirements. \QED \almostall

\prop \leproof Proof that $\forces _{P} \c(f_{\xi}, g_{\xi}) \le
	\kappa_{\xi} $:  Fix ${\xi}$.  We will write $f$ for $f_{\xi}$,
	etc. 

 Let $$A:= \bigcup \{A_\zeta:   \kappa_{\zeta}
         \le \kappa_{\xi} \}. $$  
We will show that the $g$-slaloms from $V^{P\on A}$ already
	cover~$\prod f$.  This is sufficient, because $\forces_P
	(2^{\aleph_0} )^{V^{P\on A} } \le |A| = \kappa_\xi$.

Let $p_0$ be an arbitrary condition. 
Let $\ntau$ be a name of a function~$< f$.  Find a condition $p
	\stronger p_0$ satisfying \:{finite}--\:{almostall}.

For each $l$ we now define sets $F_{\eta_l} \subseteq
\sux_{p(\alpha_l)}(\eta_l)$ as follows: 
\begindent
\ite 1 If $\alpha_l\in A$, then $F_{\eta_l}=\sux_{p(\alpha_l)}(\eta_l)$.  
\ite 2 If  $f_{\zeta_l}({k_l}) \le      g_{\xi}({k_l})/|\level_{k_l}(p)|$, then again
 $F_{\eta_l}=\sux_{p(\alpha_l)}(\eta_l)$.  
\ite 3 Otherwise, we thin out the set $\sux_{p(\alpha_l)}(\eta_l)$
such that each ${\bareta} $ in $\level_{k_l}(p) $ decides $\ntau({k_l})$ up to
at most $g({k_l})/|\level_{k_l}(p)|$ many values. 
\endent
 Here is a more detailed description of case (3): Let $k=k_l$,
	${\zeta}={\zeta}_l$.

  Note that if
neither (1) nor (2) holds, then letting $c:= f_{\xi}(k)$, $d:=
	g_{\xi}(k) / |\level_k(p)|$, we have $c/d \le h_{{\zeta}}(k)$. 

Using $(c, d)$-completeness of the norm
$\|\cdot|_{{\zeta},k}$we define a sequence
$$ \sux_{p(\alpha_l)}(\eta_l)= L(0) \supseteq L(1) \supseteq \cdots
\supseteq L(|\level_k(p)|)$$
as follows.   Let ${\bareta}_0, \ldots, {\bareta}
_{|\level_k(p)|-1}$ be an enumeration of~$\level_k(p)$. 

Given $L(i)$, we  know that for each $\nu\in L(i)$ the sequence
${\bareta}_i^{+\nu}$, (i.e., the condition 
$\trim p^{{\bareta}_i^{+\nu}} $) decides~$\ntau(k)$.  (See
\plusnotation.) 
since there only $\le c$  many possible values for $\ntau(k)$, we can
use  $(c,d)$-completeness to find a set $L(i+1) \subseteq L(i)$ and a
set $C(i)$ such that 
\begindent
\ite a $\|L(i+1)| \ge \|L(i)| - 1$
\ite b $|C(i)| \le d$.
\ite c For every $\nu\in L(i+1)$, $\trim p^{{\bareta}_i^{+\nu}} \forces
\ntau(k)\in C(i)$. 
\endent
Now let $F_{\eta_l}$ be $L(|\level_k(p)|)$, and let 
$$ ({\oplus}) \qquad\qquad B_k:= \bigcup_i C(i). $$
So $|B_k| \le |\level_k(p)|\cdot d \le g(k)$. 

Clearly $ \|F_{\eta_l}|_{{\zeta}_l, k_l} \ge \|p|_{k_l} - |\level_{k_l}(p)|
	> {1\over 2} \|p|_{k_l}$. 

This completes the definition of the sets $F_{\eta_l}$. 
\medskip

Let $q\ge p$ be the condition  defined from $p$ using the $F_{\eta_l}$ (see
\thinning).   We will find a $P\on A$-name for a $g$-slalom $\name {\bar
	B} = \<\name B_k: k\in \omega >$ such that 
$$ q \forces \hbox{$\name{\bar B} $ covers~$\ntau$.}$$

If $k$ is not a splitting level, then every ${\bareta} $  in level $k$
decides~$\ntau(k)$ by \:{tau(k)}.   So in this case we can let 
$$ B_k := \{ i: \exists {\bareta} \in \level_k(p), \trim p^{\bareta}
\forces \ntau (k)=i\}$$
This set is of size $\le |\level_k(p)| < g(k)$, and clearly $q \forces
\ntau(k)\in B_k$.

If $k$ is a splitting level, $k=k_l$, then there are three cases.

Case 1: $\alpha_l \in A$: We define $\name B_k$ to be a  $P\on A$-name
satisfying the following: 
$$ \forces_{P\on A} \name B_k = \{ i: \exists {\bareta} \in
   \level_{k+1} (p), V \thinks \trim p^{\bareta} \forces \ntau (k)=i,
   {\bareta}(\alpha_l) \subseteq r_{\alpha_l}\}$$ 
Thus, we only admit those ${\bareta} $ which agree with the generic
   real added by the forcing~$Q_{\alpha_l}$.   Clearly 
$ \forces_{P\on A} |B_k| \le \level_{k}(p) < 
   g(k)$, and $p\forces_P \ntau(k)\in B_k$.

Case 2: $f_{\zeta_l}(k) \le      g_{\xi}(k)/|\level_k(p)|$.

So we have $|\level_{k+1}(p)|\le f_{\zeta_l}(k) \cdot |\level_k(p)|
   \le g(k)$, so we can let 
$$ B_k := \{ i: \exists {\bareta} \in \level_{k+1} (p), \trim p^{\bareta}
\forces \ntau (k)=i\}$$
This set is of size $\le |\level_{k+1}(p)| \le g(k)$, and again~$p \forces
\ntau(k)\in B_k$.

Case 3:  Otherwise.   We have already defined~$B_{k_l}$ in
	$(\oplus)$.    By
	condition (c) above, $q \forces \ntau(k) \in B_k$.

\bigskip
So  indeed $q \forces$``$ \name {\bar B} = \< \name B_k:k\in \omega>$ is a
	$g$-slalom covering $\ntau$'' \hfill\hbox{\QED \leproof\
        \QED\mainTHEO\ \QED{{\rm[GSh\ 448]}}}

\bigskip

\goodbreak
\bigskip\bigskip
\bigskip\bigskip
\goodbreak

\nextinput biblio448 {\rm REFERENCES}

\def\paper#1,{{\it #1}, } 
\def\inbook#1,{in:  #1, }  \def\nextreference#1\par{\bigskip \noindent#1}
\def\journall#1:#2(#3)#4--#5.{#1, {\bf #2} (#3), pp.#4--#5.}
\def\journalp#1:#2(#3)#4.{#1, vol~#2 (#3), p.~#4.}
\def\book#1,{{\it #1},} \def\vol#1 {Vol.~#1}
\def\pages#1--#2{pp.#1--#2}
\def\by#1:{#1, }

\def\handbookST{\inbook Handbook of Set-Theoretic Topology,  ed.\ by
        K.~Kunen and J.E.~Vaughan,  North-Holland,
        Am\-ster\-dam-New York--Oxford 1984}

\def\MSRI{\inbook Proceedings of the MSRI Logic Year 1989/90, ed.\ by
	H.~Judah, W.~Just, W.~H.~Woodin}

\parindent=0cm

\font\smc=cmcsc10 

\def\key#1[#2]{\medskip\hangindent=1cm\hangafter1
{\smc [#2]}\write\mgfile{\def\string#1{#2}}}

\key\Blass[{{\sl Blass\/}}] \by A.~Blass: \paper Simple cardinal invariants,
preprint.

\key\vanDouwen[{{\sl van Douwen\/}}] \by E.K.~van Douwen: \paper The
        integers and Topology, \handbookST

\key\ComfortNegrepontis[{{\sl Comfort-Negrepontis}}] \by Comfort and
   Negrepontis: \book Theory of ultrafilters, Springer Verlag, Berlin
	Heidelberg New~York, 1974. 

\key\Miller[{{\sl Miller}}] \by A.~Miller:
\paper Some properties of measure and category, Transactions of the
AMS 266.

\key\ShelahCCCXXVI[{{\sl Shelah 326\/}}] \by S.~Shelah: \paper Vive la
difference!, to appear \MSRI. 

\key\ShelahCDXLVIIIa[{{\sl Shelah 448a}}] \by S.~Shelah: \paper Notes
on many cardinal invariants, May 1991.

\key\Vaughan[{{\sl Vaughan\/}}] \by J.E.~Vaughan: \paper Small
        uncountable cardinals and topology, \inbook Open
        problems in Topology, ed.\ by J.\ van Mill and G.\
        Reeds.

\bigskip
\bigskip
\baselineskip=12pt
\settabs10\columns
\+& 	Martin Goldstern	&&&&&Saharon Shelah\cr
\+&	Dept of Mathematics	&&&&&Dept of Mathematics\cr
\+&	Bar Ilan University	&&&&&Hebrew University\cr
\+&	52900 Ramat Gan		&&&&&Givat Ram, Jerusalem\cr
\+& 	{\tt goldstrn@bimacs.cs.biu.ac.il}
				&&&&&{\tt shelah@math.huji.ac.il}\cr

\bye